\documentclass[12pt]{amsart}

\usepackage{amssymb,mathrsfs,amsmath,amsthm,color,bm,mathtools,bbm,wasysym,cases,mathdots}
\usepackage[pagebackref]{hyperref}
\usepackage[noadjust]{cite}

\usepackage{enumerate}
\usepackage[centering]{geometry}
\newtheorem{theorem}{Theorem}

\theoremstyle{remark}
\newtheorem{remark}{Remark}

\begin{document}
	\title[Counterexamples to the inhomogeneous Duffin--Schaeffer conjecture]{Counterexamples to the inhomogeneous Duffin--Schaeffer conjecture for a residual set of shifts}
\author{Yubin He}

\address{Department of Mathematics, Shantou University, Shantou, Guangdong, 515063, China}

	\email{ybhe@stu.edu.cn}

	\author{Lingmin Liao}

	\address{School of Mathematics and Statistics, Wuhan University, Wuhan, Hubei 430072, China}

	\email{lmliao@whu.edu.cn}

\begin{abstract}
	Let $\theta\in\mathbb Q\setminus\{0\}$.  We construct an approximating
	function $\psi:\mathbb N\to[0,\frac12)$ for which
	\[
	\sum_{q=1}^{\infty}\frac{\varphi(q)}{q}\psi(q)=\infty,
	\]
	but the set of $x\in[0,1]$ for which
	\[
	|qx-a-\theta|<\psi(q),\qquad \gcd(a,q)=1,
	\]
	holds for infinitely many $(a,q)\in\mathbb Z\times\mathbb N$ has
	Lebesgue measure zero.  Thus the
	inhomogeneous analogue of the Duffin--Schaeffer conjecture fails for
	every nonzero rational shift.  By a modification of the construction, we further show that the set
	of shifts for which the inhomogeneous Duffin--Schaeffer conjecture
	fails is residual in $\mathbb R$.
\end{abstract}
	\maketitle

	\section{Introduction}
A starting point in metric Diophantine approximation is
Khintchine's theorem \cite{Khin}.  For an approximating function
\[
\psi:\mathbb N\longrightarrow[0,\infty),
\]
define
\[
W(\psi)
:=
\left\{
x\in[0,1]:
\left|x-\frac aq\right|
<
\frac{\psi(q)}{q}
\text{ for infinitely many }
(a,q)\in\mathbb Z\times\mathbb N
\right\}.
\]
Let $\lambda$ denote the Lebesgue measure on $[0,1]$.  Khintchine's
theorem \cite{Khin} states that
\[
\lambda(W(\psi))
=
\begin{cases}
	0,
	& \displaystyle\text{if }
	\sum_{q=1}^{\infty}\psi(q)<\infty,\\[2mm]
	1,
	& \displaystyle\text{if }
	\sum_{q=1}^{\infty}\psi(q)=\infty
	\text{ and $\psi$ is non-increasing}.
\end{cases}
\]
The convergence case does not require the monotonicity assumption.
It is therefore natural to ask whether monotonicity can also be removed
in the divergence case.  Duffin and Schaeffer \cite{DS} showed that the monotonicity assumption
in Khintchine's theorem cannot simply be removed.  More precisely, they
constructed an approximating function $\psi$ for which
\[
\sum_{q=1}^{\infty}\psi(q)=\infty,
\]
but nevertheless
\[
\lambda(W(\psi))=0.
\]
This led Duffin and Schaeffer to restrict attention to reduced
fractions $a/q$, with $(a,q)=1$, and to replace by
\begin{equation}\label{eq:sum}
	\sum_{q=1}^{\infty}\frac{\varphi(q)}{q}\psi(q)=\infty.
\end{equation}
They conjectured that the divergence of the sum \eqref{eq:sum} is sufficient to
ensure that, for $\lambda$-almost every $x\in[0,1]$,
\[
\left|x-\frac aq\right|<\frac{\psi(q)}q,
\qquad \gcd(a,q)=1,
\]
has infinitely many solutions.

The Duffin--Schaeffer conjecture was proved by Koukoulopoulos and
Maynard \cite{KM}.  Thus, for reduced rational approximation, no
monotonicity assumption on $\psi$ is required.

A natural inhomogeneous analogue is obtained by introducing a fixed
shift $\theta\in\mathbb R$ and replacing the rational centres $a/q$
by
\[
\frac{a+\theta}{q},
\]
while retaining the reducedness condition $(a,q)=1$.  Accordingly,
define
\[
W_\theta^*(\psi)
:=
\left\{
x\in[0,1]:
\begin{aligned}
	& \bigg|x-\frac{a+\theta}{q}\bigg|<\frac{\psi(q)}{q}
	\text{ for infinitely many} \\
	& (a,q)\in\mathbb Z\times\mathbb N
	\text{ with }\gcd(a,q)=1
\end{aligned}
\right\}.
\]
The corresponding inhomogeneous Duffin--Schaeffer conjecture is as follows.

\medskip
\noindent\textbf{Inhomogeneous Duffin--Schaeffer conjecture}.
\emph{Let $\theta\in\mathbb R$ and let
	$\psi:\mathbb N\to[0,\infty)$. Then,
	\[
	\lambda\bigl(W_\theta^*(\psi)\bigr)
	=
	\begin{cases}
		0,
		& \displaystyle\text{if }
		\sum_{q=1}^{\infty}
		\frac{\varphi(q)}{q}\psi(q)<\infty,\\[3mm]
		1,
		& \displaystyle\text{if }
		\sum_{q=1}^{\infty}
		\frac{\varphi(q)}{q}\psi(q)=\infty.
	\end{cases}
	\]
}

Several related forms of the inhomogeneous Duffin--Schaeffer conjecture
have been studied.  For a rational shift
\[
\theta=\frac{A}{B},
\qquad \gcd(A,B)=1,
\]
Beresnevich, Hauke and Velani \cite{BHV} considered the set of
$x\in[0,1]$ for which
\[
\left|x-\frac{a+A/B}{q}\right|<\frac{\psi(q)}{q},
\qquad \gcd(A+aB,q)=1,
\]
holds for infinitely many $(a,q)\in\mathbb Z\times\mathbb N$.  They
proved that this set has full Lebesgue measure whenever
\[
\sum_{q=1}^{\infty}\frac{\varphi(q)}{q}\psi(q)=\infty.
\]
They also established the weak inhomogeneous Duffin--Schaeffer
conjecture for every rational shift, where no coprimality condition is
imposed.  These results are related to, but distinct from, the
inhomogeneous Duffin--Schaeffer conjecture considered here, in which
the coprimality condition is $\gcd(a,q)=1$.

A different problem arises when the inhomogeneous shift is allowed
to vary with the denominator.  More precisely, Hauke and Ram\'irez
\cite{HR} considered a sequence of shifts
\[
\mathbf y=(y_q)_{q\in\mathbb N}
\]
and the corresponding inequalities
\begin{equation}\label{eq:moving}
	|qx-a-y_q|<\psi(q),
	\qquad (a,q)=1.
\end{equation}
They constructed a sequence $\mathbf y$ and an approximating function
$\psi$ such that
\[
\sum_{q=1}^{\infty}\frac{\varphi(q)}{q}\psi(q)=\infty,
\]
while the set of $x$ for which the inequality \eqref{eq:moving} has infinitely
many solutions has Lebesgue measure zero.  Thus the one-dimensional
Duffin--Schaeffer statement fails when the inhomogeneous shift is
allowed to depend on $q$.  The construction of Hauke and Ram\'irez, however, does not
produce a counterexample when the shift is fixed, that is, when
$y_q=\theta$ for every $q$. Further results directly concerning the inhomogeneous Duffin--Schaeffer conjecture can be found in \cite{CT,Ramirez,Yu}.

The first result of this paper shows that the inhomogeneous Duffin--Schaeffer conjecture fails for every nonzero rational shift.

\begin{theorem}\label{t:main}
	Let
	\[
	\theta=\frac AB\in\mathbb Q\setminus\{0\}
	\qquad\text{with\quad}
	A\in\mathbb Z\setminus\{0\},\quad
	B\in\mathbb N,\quad
	\gcd(A,B)=1.
	\]
	Then there exists a function
	\[
	\psi:\mathbb N\to[0,\tfrac12)
	\]
	such that
	\[
	\sum_{q=1}^{\infty}
	\frac{\varphi(q)}q\psi(q)=\infty,
	\]
	whereas
	\[
	\lambda(W_\theta^*(\psi))=0.
	\]
\end{theorem}

\begin{remark}
	The reason for the failure in the inhomogeneous setting can already be
	seen from a difference with the homogeneous case.
	 In the homogeneous setting, the coprimality condition $(a,q)=1$
	 eliminates repetitions arising from different representations of the
	 same rational number.  Indeed, if
	 \[
	 \frac{a_1}{q_1}=\frac{a_2}{q_2},
	 \qquad
	 \gcd(a_1,q_1)=\gcd(a_2,q_2)=1,
	 \]
	 then necessarily
	 \[
	 a_1=a_2
	 \qquad\text{and}\qquad
	 q_1=q_2.
	 \]
	 Thus two distinct coprime pairs cannot represent the same rational
	 number.

	In the inhomogeneous setting, this property is lost.  Even when
	$\gcd(a,q)=1$, distinct pairs $(a,q)$ may give the same value of
	\[
	\frac{a+\theta}{q}.
	\]
	For example, take $\theta=1$.  For every $q\ge 2$, let
	\[
	a=q-1.
	\]
	Then,
	\[
	\gcd(a,q)=\gcd(q-1,q)=1,
	\]
	but
	\[
	\frac{a+\theta}{q}
	=
	\frac{q-1+1}{q}
	=
	1.
	\]
	Thus the infinitely many distinct coprime pairs
	\[
	(q-1,q),\qquad q\ge2,
	\]
	all give the same value.
	The construction of the counterexample in Theorem~\ref{t:main} is
	based precisely on this phenomenon: many distinct coprime pairs $\gcd(a,q)=1$ may
	produce the same rational $(a+\theta)/q$.
\end{remark}
\begin{remark}
	In higher dimensions, Pollington and Vaughan \cite{PV} proved the homogeneous Duffin--Schaeffer conjecture, while Hauke and Ram\'irez \cite{HR} established the corresponding inhomogeneous result in dimensions $m\ge3$. By contrast, Theorem~\ref{t:main} shows that in dimension one the fixed-shift conjecture fails for every nonzero rational shift.
\end{remark}

A suitable modification of the proof of Theorem~\ref{t:main} allows us
to extend the result considerably. Recall that a subset of $\mathbb R$
is called \emph{residual} if its complement is of first category, that
is, if its complement is a countable union of nowhere dense sets.
\begin{theorem}\label{t:residual}
	The set of shifts $\theta\in\mathbb R$ for which the inhomogeneous Duffin--Schaeffer conjecture fails
	is residual in $\mathbb R$.
\end{theorem}

	\section{Proof of Theorem \ref{t:main}}

	Fix
	\[
	\theta=\frac AB\neq 0
	\qquad\text{with}\quad
	A\in\mathbb Z\setminus\{0\},\quad
	B\in\mathbb N,\quad
	\gcd(A,B)=1.
	\]
	Let
	\[
	\ell_1<\ell_2<\ell_3<\cdots
	\]
	be the increasing enumeration of all odd primes not dividing $AB$, and set
	\[
	\mathcal P_m:=\{\ell_1,\ldots,\ell_m\},
	\qquad
	P_m:=\prod_{p\in\mathcal P_m}p
	=\ell_1\cdots\ell_m.
	\]
	Thus $P_m$ is squarefree. Since the sum of the reciprocals of all
	primes diverges, and only finitely many primes are excluded, namely
	$2$ and the prime divisors of $AB$, we still have
	\[
	\sum_{m=1}^{\infty}\frac1{\ell_m}=\infty.
	\]

	For each $m\ge 1$, define the denominator block
	\[
	\mathcal Q_m
	:=
	\{q:q\mid P_m \text{ and } \ell_m\mid q\}.
	\]
	These blocks are pairwise disjoint. Indeed, if $m<n$, then every
	$q\in\mathcal Q_m$ has all of its prime divisors among
	$\ell_1,\ldots,\ell_m$, whereas every $q\in\mathcal Q_n$ is divisible
	by the new prime $\ell_n$.

	For any $q\in\mathcal Q_m$, define
	\[
	\psi(q):=\frac{q}{4P_m}.
	\]
	Then,
	\begin{equation}\label{eq:sum block}
		\sum_{q\in\mathcal Q_m}\frac{\varphi(q)}q\psi(q)
		=
		\frac{1}{4P_m}
		\sum_{\substack{q\mid P_m\\ \ell_m\mid q}}\varphi(q).
	\end{equation}
	Note that every divisor $q\mid P_m$ with $\ell_m\mid q$ can be written uniquely as
	\[
	q=\ell_m e,
	\qquad \text{with}\quad e\mid P_{m-1}\quad\text{and}\quad\gcd(\ell_m,e)=1.
	\]
	Then, the
	multiplicativity of Euler's totient function gives
	\[
	\varphi(q)
	=
	\varphi(\ell_m e)
	=
	\varphi(\ell_m)\varphi(e)
	=
	(\ell_m-1)\varphi(e).
	\]
	Therefore, by the classical identity
	$\sum_{d\mid n}\varphi(d)=n$,
	\begin{align}
		\sum_{q\in\mathcal Q_m}\frac{\varphi(q)}q\psi(q)
		&=
		\frac{1}{4P_m}
		(\ell_m-1)
		\sum_{e\mid P_{m-1}}\varphi(e) \notag\\
		&=
		\frac{1}{4P_m}
		(\ell_m-1)P_{m-1}.\notag
	\end{align}
	By the relation $P_m=\ell_mP_{m-1}$,
	\begin{equation}\label{eq:block-mass-single}
		\frac{1}{4P_m}
		(\ell_m-1)P_{m-1}
		=
		\frac{\ell_m-1}{4\ell_m}
		=
		\frac {1}{4}
		\left(1-\frac1{\ell_m}\right).
	\end{equation}
	 Since $\ell_m\ge 3$, it follows that
	\[
	1-\frac1{\ell_m}\ge\frac23,
	\]
	and hence
\begin{equation}\label{eq:computation}
	\sum_{q\in\mathcal Q_m}\frac{\varphi(q)}q\psi(q)
	\ge
	\frac{1}{6}.
\end{equation}
	Thus every block contributes a uniformly positive amount to the
	sum \eqref{eq:sum block}.

For each denominator $q$ and each integer $a$ with $\gcd(a,q)=1$, define
the corresponding \emph{approximation centre} by
\[
\frac{a+\theta}{q}
=
\frac{a+A/B}{q}.
\]
The associated \emph{approximation interval} is
\[
I_{a,q}
=
\left\{
x\in[0,1]:
\bigg|x-\frac{a+\theta}{q}\bigg|<\frac{\psi(q)}{q}
\right\}.
\]
We next estimate the total measure of the approximation intervals contributed
by each block $\mathcal Q_m$.
 For any $q\in\mathcal Q_m$,
since $q\mid P_m$, there is a unique divisor $d\mid P_m$ such that
\[
q=\frac{P_m}{d}.
\]
An approximation centre is therefore of the form
\[
\frac{a+A/B}{q}
=
\frac{d(Ba+A)}{BP_m}.
\]
Set
\[
n:=d(Ba+A).
\]
Then, all approximation centres arising from $\mathcal Q_m$ lie on the
common lattice
\[
\frac{1}{BP_m}\mathbb Z\Big/\mathbb Z.
\]
Accordingly, we refer to $n=d(Ba+A)$ as a \emph{lattice numerator}
associated with the block $\mathcal Q_m$.

We now examine the congruence restrictions imposed on $n$ by the
condition $\gcd(a,q)=1$.  Let $p\in\mathcal P_m$.  Since $P_m$ is
squarefree, exactly one of
\[
p\mid d
\qquad\text{or}\qquad
p\mid \frac{P_m}{d}=q
\]
holds.

If $p\mid d$, then
\[
n=d(Ba+A)\equiv0\pmod p.
\]
If $p\mid P_m/d=q$, then $p\nmid d$. Since $\gcd(a,q)=1$, we have $p\nmid a$.  Moreover, by the choice of the primes $\ell_i$,
$p\nmid B$.  Hence,
\[
n-Ad=dBa\not\equiv0\pmod p,
\]
and therefore,
\[
n\not\equiv Ad\pmod p.
\]
Thus, every lattice numerator $n$ associated with the block $\mathcal Q_m$ satisfies
\[
\begin{cases}
	n\equiv0\pmod p, & p\mid d,\\[1mm]
	n\not\equiv Ad\pmod p, & p\mid P_m/d.
\end{cases}
\]

For the next counting argument, we retain only the congruence information
modulo the prime divisors of $P_m$.  Let $r$ denote the residue class
of $n$ modulo $P_m$; that is,
\[
r\in\mathbb Z/P_m\mathbb Z\quad\text{and}\quad r=n \pmod{P_m}.
\]
Since every $p\in\mathcal P_m$ divides $P_m$, the congruence
$r\equiv n\pmod{P_m}$ implies
\[
r\equiv n\pmod p
\qquad\text{for every }p\in\mathcal P_m.
\]
Therefore, the congruence restrictions satisfied by $n$ pass directly
to $r$.  In particular,
\begin{equation}\label{eq:witness}
	\begin{cases}
		r\equiv0\pmod p, & p\mid d,\\[1mm]
		r\not\equiv Ad\pmod p, & p\mid P_m/d.
	\end{cases}
\end{equation}
Whenever a divisor $d\mid P_m$ satisfies \eqref{eq:witness} for a
residue class $r\pmod{P_m}$, we call $d$ a \emph{witness} for $r$.

Define
\[
R_A(P_m)
:=
\left\{
r\in\mathbb Z/P_m\mathbb Z:
\begin{array}{l}
	\text{there exists }d\mid P_m\text{ such that}\\
	r\equiv0\pmod p
	\quad\text{for every }p\mid d,\\
	r\not\equiv Ad\pmod p
	\quad\text{for every }p\mid P_m/d
\end{array}
\right\}.
\]
Thus, every lattice numerator $n$ associated with the block $\mathcal Q_m$ determines a residue class $r\pmod{P_m}$ belonging
to $R_A(P_m)$. Notice, however, that the converse need not hold.  Indeed, if
$q\in\mathcal Q_m$ and
\[
q=\frac{P_m}{d},
\]
then the defining condition $\ell_m\mid q$ forces
\[
\ell_m\nmid d.
\]
In the definition of $R_A(P_m)$, by contrast, we allow every divisor
$d\mid P_m$, including those divisible by $\ell_m$.  Hence, $R_A(P_m)$ may contain residue classes that do not correspond to
any lattice numerator $n$ associated with the block $\mathcal Q_m$.
This enlargement is convenient and
causes no difficulty, since we only use $R_A(P_m)$ to obtain an upper
bound for the number of possible approximation centres.

Define
\[
\rho_A(P_m):=\frac{|R_A(P_m)|}{P_m},
\]
where $|\cdot|$ denotes the cardinality of a finite set.
Since there are exactly $P_m$ residue classes modulo $P_m$,
$\rho_A(P_m)$ measures the proportion of residue classes satisfying
the congruence restrictions above.

Note that each residue class $r\pmod{P_m}$ has exactly
$B$ lifts modulo $BP_m$:
\[
r,\ r+P_m,\ \ldots,\ r+(B-1)P_m
\pmod{BP_m}.
\]
Hence, the number of lattice numerators
associated with the block $\mathcal Q_m$ is at most
\[
B|R_A(P_m)|.
\]
Since each approximation interval has radius $1/(4P_m)$ and therefore
length $1/(2P_m)$, we obtain
\begin{equation}
	\lambda\bigg(\bigcup_{q\in\mathcal Q_m}\bigcup_{\gcd(a,q)=1}B\bigg(\frac{a+A/B}{q},\frac{1}{4P_m}\bigg)\bigg)
	\le
	B|R_A(P_m)|\frac{1}{2P_m}
	=B\rho_A(P_m)/2.
	\label{eq:block-measure-single}
\end{equation}

It remains to prove that
\[
\rho_A(P_m)\longrightarrow0.
\]
For each residue class
\[
r\in\mathbb Z/P_m\mathbb Z,
\]
define its \emph{zero-pattern} by
\[
Z(r)
:=
\{p\in\mathcal P_m:r\equiv0\pmod p\}.
\]
Thus $Z(r)$ records precisely the primes in $\mathcal P_m$ modulo which $r$ vanishes.

For each subset $Z\subseteq\mathcal P_m$, including $Z=\emptyset$, define
\[
C_Z
:=
\{r\in\mathbb Z/P_m\mathbb Z:Z(r)=Z\}.
\]
Equivalently,
\[
C_Z
=
\left\{
r\pmod{P_m}:
\begin{cases}
	r\equiv0\pmod p, & p\in Z,\\
	r\not\equiv0\pmod p, & p\in\mathcal P_m\setminus Z
\end{cases}
\right\}.
\]
Every residue class modulo $P_m$ has a unique zero-pattern.  Hence the
sets $C_Z$ are pairwise disjoint and form a partition:
\begin{equation}\label{eq:disjoint}
	\mathbb Z/P_m\mathbb Z
	=
	\bigsqcup_{Z\subseteq\mathcal P_m}C_Z.
\end{equation}
Therefore, intersecting this disjoint decomposition with $R_A(P_m)$ gives
\[
R_A(P_m)
=
\bigsqcup_{Z\subseteq\mathcal P_m}
\bigl(R_A(P_m)\cap C_Z\bigr),
\]
and consequently
\begin{equation}\label{eq:estimate}
	|R_A(P_m)|
	=
	\sum_{Z\subseteq\mathcal P_m}
	|R_A(P_m)\cap C_Z|.
\end{equation}

Let $Z\subseteq \mathcal P_m$. We shall derive two upper bounds for $|R_A(P_m)\cap C_Z|$.  The first follows directly
from the size of $C_Z$, while the second uses the restrictions imposed
by the possible witnesses.

We begin with the first bound by estimating the size of $C_Z$. Since
\[
P_m=\prod_{p\in\mathcal P_m}p
\]
is squarefree, its prime factors are pairwise coprime.  Hence the
Chinese remainder theorem gives a bijection
\[
\Phi:\mathbb Z/P_m\mathbb Z
\longrightarrow
\prod_{p\in\mathcal P_m}\mathbb Z/p\mathbb Z
\]
defined by
\begin{equation}\label{eq:CRT}
	\Phi\bigl(r\pmod{P_m}\bigr)
	=
	\bigl(r\pmod p\bigr)_{p\in\mathcal P_m}.
\end{equation}
In other words, a residue class modulo $P_m$ is uniquely determined
by its residues modulo the primes $p\in\mathcal P_m$.

If $p\in Z$, then every $r\in C_Z$ satisfies
\[
r\equiv0\pmod p,
\]
so there is exactly one possible residue modulo $p$.  If
$p\in\mathcal P_m\setminus Z$, then
\[
r\not\equiv0\pmod p,
\]
so there are exactly $p-1$ possibilities modulo $p$.  Hence,
\begin{equation}\label{eq:upper 1}
	|C_Z|
	=
	\prod_{p\in Z}1
	\prod_{p\in\mathcal P_m\setminus Z}(p-1)
	=
	\prod_{p\in\mathcal P_m\setminus Z}(p-1).
\end{equation}

We next derive a second bound for $|R_A(P_m)\cap C_Z|$ by exploiting the restrictions
on the possible witnesses.
Let
\[
r\in R_A(P_m)\cap C_Z,
\]
and let $d\mid P_m$ be a witness for $r$.  Then
\[
r\equiv0\pmod p
\qquad\text{for every }p\mid d.
\]
Since $r\in C_Z$,
\[
r\equiv0\pmod p
\quad\Longleftrightarrow\quad
p\in Z.
\]
It follows that every prime divisor of $d$ belongs to $Z$.

Define
\[
D_Z:=\prod_{p\in Z}p,
\]
with the convention $D_\varnothing=1$.  Since $d\mid P_m$ and $P_m$
is squarefree, we obtain
\[
d\mid D_Z.
\]
Thus every possible witness $d$ is a divisor of $D_Z$.  Since $D_Z$
is squarefree with $|Z|$ prime factors, it has exactly
\begin{equation}\label{eq:divisors DZ}
	2^{|Z|}
\end{equation}
positive divisors.  Hence, for a fixed $Z$, there are at most
$2^{|Z|}$ possible witnesses.

Now fix a divisor
\[
d\mid D_Z
\]
and define
\[
C_Z(d)
:=
\{r\in C_Z:\text{$d$ is a witness for $r$}\}.
\]
The Chinese remainder theorem \eqref{eq:CRT} allows us to determine the size of $C_Z(d)$ by considering the residue of $r$
modulo each $p\in\mathcal P_m$.

If $p\mid d$, then $p\in Z$, and hence
\[
r\equiv0\pmod p.
\]
Thus there is exactly one possibility modulo $p$.

If $p\in Z$ but $p\nmid d$, then again
\[
r\equiv0\pmod p.
\]
Since $p\nmid A$ and $p\nmid d$,
\[
Ad\not\equiv0\pmod p,
\]
so the condition
\[
r\not\equiv Ad\pmod p
\]
is automatically satisfied.  Thus there is again exactly one
possibility modulo $p$.

Finally, suppose that
\[
p\in\mathcal P_m\setminus Z.
\]
Since $d\mid D_Z$, we have $p\nmid d$.  The condition $r\in C_Z$
requires
\[
r\not\equiv0\pmod p,
\]
while the condition that $d$ is a witness requires
\[
r\not\equiv Ad\pmod p.
\]
The residues $0$ and $Ad$ are distinct modulo $p$, since
$p\nmid A$ and $p\nmid d$.  Therefore exactly two residue classes are
excluded, leaving
\[
p-2
\]
possibilities modulo $p$.

Since $P_m$ is squarefree, the Chinese remainder theorem \eqref{eq:CRT} combines
these choices uniquely into a residue class modulo $P_m$.  Therefore
\begin{equation}\label{eq:upper 2}
	|C_Z(d)|
	=
	\prod_{p\in Z}1
	\prod_{p\in\mathcal P_m\setminus Z}(p-2)
	=
	\prod_{p\in\mathcal P_m\setminus Z}(p-2).
\end{equation}

We now obtain two upper bounds for $|R_A(P_m)\cap C_Z|$.  First,
\[
|R_A(P_m)\cap C_Z|
\le
|C_Z|
\stackrel{\eqref{eq:upper 1}}{=}
\prod_{p\in\mathcal P_m\setminus Z}(p-1).
\]
Second, there are at most $2^{|Z|}$ possible witnesses (see \eqref{eq:divisors DZ}), and for each
fixed witness $d$,
\[
|C_Z(d)|
\stackrel{\eqref{eq:upper 2}}{=}
\prod_{p\in\mathcal P_m\setminus Z}(p-2).
\]
Hence
\[
|R_A(P_m)\cap C_Z|
\le
2^{|Z|}
\prod_{p\in\mathcal P_m\setminus Z}(p-2).
\]
Combining the two estimates,
\[
|R_A(P_m)\cap C_Z|
\le
\min\left\{
\prod_{p\in\mathcal P_m\setminus Z}(p-1),
\;
2^{|Z|}
\prod_{p\in\mathcal P_m\setminus Z}(p-2)
\right\}.
\]
Using
\[
\min\{x,y\}\le\sqrt{xy},
\]
we obtain
\[
|R_A(P_m)\cap C_Z|
\le
2^{|Z|/2}
\prod_{p\in\mathcal P_m\setminus Z}
\sqrt{(p-1)(p-2)}.
\]
Substitution of this upper bound into \eqref{eq:estimate} gives
\begin{align}
	|R_A(P_m)|
	&\le
	\sum_{Z\subseteq\mathcal P_m}
	2^{|Z|/2}
	\prod_{p\in\mathcal P_m\setminus Z}
	\sqrt{(p-1)(p-2)}. \notag\\
	&=\sum_{Z\subseteq\mathcal P_m}
	\prod_{p\in Z}\sqrt{2}
	\prod_{p\in\mathcal P_m\setminus Z}
	\sqrt{(p-1)(p-2)}.\notag
\end{align}
To estimate the last sum, observe that for each prime
$p\in\mathcal P_m$ there are two possibilities.  If $p\in Z$, then
the corresponding factor is $\sqrt2$, whereas if
$p\in\mathcal P_m\setminus Z$, the corresponding factor is $\sqrt{(p-1)(p-2)}$.
Therefore, summing over all subsets $Z\subseteq\mathcal P_m$ is
exactly the expansion of the product
\[
\prod_{p\in\mathcal P_m}
\left(
\sqrt2+\sqrt{(p-1)(p-2)}
\right).
\]
Therefore,
\[	|R_A(P_m)|\le\prod_{p\in\mathcal P_m}
\left(
\sqrt2+\sqrt{(p-1)(p-2)}
\right).\]
Consequently,
\begin{equation}
	\rho_A(P_m)=\frac{|R_A(P_m)|}{P_m}
	\le
	\prod_{p\in\mathcal P_m}
	\frac{\sqrt2+\sqrt{(p-1)(p-2)}}{p}.
	\label{eq:rho-product-single}
\end{equation}

For a prime $p$, set
\[
F(p)
:=
\frac{\sqrt2+\sqrt{(p-1)(p-2)}}{p}.
\]
Since
\[
\frac{\sqrt{(p-1)(p-2)}}{p}
=
\sqrt{1-\frac3p+\frac2{p^2}},
\]
the elementary inequality
\[
\sqrt{1-x}\le1-\frac{x}{2}
\qquad (0\le x\le1)
\]
gives
\[
F(p)
\le
1-\frac{\delta}{p}+\frac1{p^2},
\qquad
\delta:=\frac32-\sqrt2>0.
\]
For all sufficiently large $p$,
\[
\frac1{p^2}\le\frac{\delta}{2p},
\]
and hence
\[
F(p)
\le
1-\frac{\delta}{2p}
\le
\exp\left(-\frac{\delta}{2p}\right).
\]
It follows that there exists a constant $c>0$ such that
\[
\rho_A(P_m)
\le
c
\exp\bigg(
-\frac{\delta}{2}
\sum_{\substack{1\le i\le m}}
\frac1{\ell_i}
\bigg).
\]
Since
\[
\sum_{i=1}^{\infty}\frac1{\ell_i}=\infty,
\]
the exponent tends to $-\infty$, and therefore
\begin{equation}\label{eq:rho0}
	\rho_A(P_m)\longrightarrow0.
\end{equation}

We may therefore choose a subsequence
\[
m_1<m_2<m_3<\cdots
\]
such that
\[
\rho_A(P_{m_j})\le2^{-j}
\qquad (j\ge1).
\]
We now define the approximating function
\[
\psi(q)
=
\begin{cases}
	\displaystyle\frac{q}{4P_{m_j}},
	&
	q\in\mathcal Q_{m_j}\text{ for some }j,\\[3mm]
	0,
	&
	\text{otherwise}.
\end{cases}
\]
Since the blocks $\mathcal Q_{m_j}$ are pairwise disjoint, this definition is
unambiguous.  Moreover, $q\le P_{m_j}$ implies
\[
\psi(q)\le\frac {1}{4}<\frac12.
\]
By \eqref{eq:block-mass-single},
\begin{align*}
	\sum_{q=1}^{\infty}\frac{\varphi(q)}q\psi(q)
	&=
	\sum_{j=1}^{\infty}
	\sum_{q\in\mathcal Q_{m_j}}
	\frac{\varphi(q)}q\psi(q)=
	\sum_{j=1}^{\infty}
	\frac {1}{4}
	\left(1-\frac1{\ell_{m_j}}\right)\\
	&\ge
	\sum_{j=1}^{\infty}\frac{1}{6}
	=
	\infty.
\end{align*}

On the other hand, \eqref{eq:block-measure-single} gives
\[
\sum_{j=1}^{\infty}\lambda\bigg(\bigcup_{q\in\mathcal Q_{m_j}}\bigcup_{\gcd(a,q)=1}B\bigg(\frac{a+A/B}{q},\frac{1}{4P_{m_j}}\bigg)\bigg)
\le
\sum_{j=1}^{\infty}2^{-(j+1)}B
<
\infty.
\]
By the Borel--Cantelli lemma,
\[
\lambda(W_\theta^*(\psi))=0.
\]
This proves Theorem~\ref{t:main}.

\begin{remark}
	The assumption $A\neq0$ is essential in the counting argument above.
	Indeed, suppose that $p\in Z$ but $p\nmid d$.  Since $r\in C_Z$, we have
	\[
	r\equiv0\pmod p.
	\]
	At the same time, because $d$ is a witness for $r$, we require
	\[
	r\not\equiv Ad\pmod p.
	\]
	When $A\neq0$, our choice of the primes in $\mathcal P_m$ ensures that
	$p\nmid A$, and since $p\nmid d$ we have
	\[
	Ad\not\equiv0\pmod p.
	\]
	Thus the two conditions are compatible.

	If $A=0$, however, the second condition becomes
	\[
	r\not\equiv0\pmod p,
	\]
	which contradicts $r\equiv0\pmod p$.  Hence the counting argument used
	for nonzero rational shifts does not extend to the homogeneous case.
	There is therefore no contradiction with the result
	of Koukoulopoulos and Maynard \cite{KM}.
\end{remark}

\section{Proof of Theorem \ref{t:residual}}
In this section, we modify the construction from the previous section
to prove Theorem~\ref{t:residual}.

For each reduced nonzero rational number $A/B$, let
\[
\ell_m(A,B),\qquad P_m(A,B),\qquad \mathcal Q_m(A,B)
\]
denote the quantities $\ell_m$, $P_m$ and $\mathcal Q_m$, respectively,
constructed in the proof of Theorem~\ref{t:main} for the shift $A/B$.
Thus, $\ell_1(A,B)<\ell_2(A,B)<\cdots$ are the odd primes not
dividing $AB$,
\[
P_m(A,B)=\prod_{i=1}^m\ell_i(A,B),
\]
and
\[
\mathcal Q_m(A,B)
=
\left\{
q:q\mid P_m(A,B)
\text{ and }\ell_m(A,B)\mid q
\right\}.
\]
We also write $\rho_A(P_m(A,B))$ for the corresponding quantity
defined in the proof of Theorem~\ref{t:main}. For every fixed
$A/B$, we have
\begin{equation}\label{eq:rho to 0}
	\rho_A(P_m(A,B))\longrightarrow0
	\qquad\text{as }m\to\infty.
\end{equation}
We first record an observation that will be used in the construction.
Fix a reduced nonzero rational number $A/B$. By \eqref{eq:rho to 0},
we can choose $m=m(A,B)$ sufficiently large so that
\begin{equation}\label{eq:residual-choice}
	\ell_{m(A,B)}(A,B)>B
	\qquad\text{and}\qquad
	\frac B2
	\rho_A\bigl(P_{m(A,B)}(A,B)\bigr)
	<
	2^{-B-1}.
\end{equation}
For simplicity, set
\[
\ell(A,B):=\ell_{m(A,B)}(A,B),\qquad
P(A,B):=P_{m(A,B)}(A,B),
\]
and
\[
\mathcal Q(A,B):=\mathcal Q_{m(A,B)}(A,B).
\]

Let
\[
\gamma\in
B\left(
\frac AB,
\frac{\ell(A,B)}{8P(A,B)}
\right).
\]
For any $q\in\mathcal Q(A,B)$, we have
$q\ge\ell(A,B)$, and hence
\[
\left|\gamma-\frac AB\right|
<
\frac{\ell(A,B)}{8P(A,B)}
\le
\frac{q}{8P(A,B)}.
\]
Therefore, if $a\in\mathbb Z$ with $\gcd(a,q)=1$ and
\[
\left|x-\frac{a+\gamma}{q}\right|
<
\frac{1}{8P(A,B)},
\]
then
\begin{align*}
	\left|x-\frac{a+A/B}{q}\right|
	&\le
	\left|x-\frac{a+\gamma}{q}\right|
	+
	\left|\frac{\gamma-A/B}{q}\right|\\
	&<
	\frac{1}{8P(A,B)}
	+
	\frac{1}{8P(A,B)}\\
	&=
	\frac{1}{4P(A,B)}.
\end{align*}
Consequently,
\begin{align}
	&\lambda\left(
	\bigcup_{q\in\mathcal Q(A,B)}
	\bigcup_{\gcd(a,q)=1}
	\left\{
	x\in[0,1]:
	\left|x-\frac{a+\gamma}{q}\right|
	<
	\frac{1}{8P(A,B)}
	\right\}
	\right)\notag\\
	\le{}&
	\lambda\left(
	\bigcup_{q\in\mathcal Q(A,B)}
	\bigcup_{\gcd(a,q)=1}
	\left\{
	x\in[0,1]:
	\left|x-\frac{a+A/B}{q}\right|
	<
	\frac{1}{4P(A,B)}
	\right\}
	\right)\notag\\
	\le{}&
	\frac B2\rho_A(P(A,B))
	<
	2^{-B-1}.
	\label{eq:residual-block}
\end{align}

For each $n\in\mathbb N$, define
\[
\mathcal G_n
:=
\bigcup_{\substack{
		A/B\in\mathbb Q\setminus\{0\}\\
		\gcd(A,B)=1,\ B\ge n}}
B\left(
\frac AB,
\frac{\ell(A,B)}{8P(A,B)}
\right).
\]
Clearly, $\mathcal G_n$ is open. Moreover, every reduced nonzero
rational number $A/B$ with $B\ge n$ belongs to $\mathcal G_n$.
Since such rational numbers are dense in $\mathbb R$, so is
$\mathcal G_n$. Therefore, by the Baire
category theorem,
\[
\mathcal G
:=
\bigcap_{n=1}^{\infty}\mathcal G_n
\]
is a dense $G_\delta$ subset of $\mathbb R$.

We now show that the inhomogeneous Duffin--Schaeffer conjecture fails
for every $\theta\in\mathcal G$.

Fix $\theta\in\mathcal G$. Since $\theta\in\mathcal G_n$ for every
$n\in\mathbb N$, we may choose a sequence of reduced nonzero rational
numbers
\[
\frac{A_j}{B_j},
\qquad j\ge1,
\]
such that
\[
\theta\in
B\left(
\frac{A_j}{B_j},
\frac{\ell(A_j,B_j)}{8P(A_j,B_j)}
\right)
\]
and
\begin{equation}\label{eq:successive-blocks}
	B_{j+1}>P(A_j,B_j)
	\qquad (j\ge1).
\end{equation}
In particular,
\[
B_1<B_2<\cdots\longrightarrow\infty.
\]
Moreover, since $\ell(A,B)>B$ by our choice of $m(A,B)$,
\eqref{eq:successive-blocks} implies
\[
\ell(A_{j+1},B_{j+1})
>
B_{j+1}
>
P(A_j,B_j).
\]
It follows that the corresponding blocks
\[
\mathcal Q(A_j,B_j),
\qquad j\ge1,
\]
are pairwise disjoint. Indeed, every $q\in\mathcal Q(A_j,B_j)$
satisfies
\[
q\le P(A_j,B_j),
\]
whereas every $q\in\mathcal Q(A_{j+1},B_{j+1})$ satisfies
\[
q\ge\ell(A_{j+1},B_{j+1})
>
P(A_j,B_j).
\]

We now define
\[
\psi(q)
=
\begin{cases}
	\displaystyle
	\frac{q}{8P(A_j,B_j)},
	&
	q\in\mathcal Q(A_j,B_j)
	\text{ for some }j,\\[3mm]
	0,
	&
	\text{otherwise}.
\end{cases}
\]
Since the blocks $\mathcal Q(A_j,B_j)$ are pairwise disjoint, this
definition is unambiguous. Moreover, for
$q\in\mathcal Q(A_j,B_j)$, we have $q\le P(A_j,B_j)$, and hence
\[
0\le\psi(q)\le\frac18<\frac12.
\]

For each $j$, the same computation as in \eqref{eq:computation} gives
\begin{align*}
	\sum_{q\in\mathcal Q(A_j,B_j)}
	\frac{\varphi(q)}q\psi(q)
	&=
	\frac{1}{8P(A_j,B_j)}
	\sum_{q\in\mathcal Q(A_j,B_j)}\varphi(q)\\
	&=
	\frac18
	\left(
	1-\frac1{\ell(A_j,B_j)}
	\right)\\
	&\ge
	\frac1{12}.
\end{align*}
Therefore,
\[
\sum_{q=1}^{\infty}
\frac{\varphi(q)}q\psi(q)
\ge
\sum_{j=1}^{\infty}\frac1{12}
=
\infty.
\]

On the other hand, since
\[
\theta\in
B\left(
\frac{A_j}{B_j},
\frac{\ell(A_j,B_j)}{8P(A_j,B_j)}
\right),
\]
the estimate \eqref{eq:residual-block}, applied with
$A=A_j$ and $B=B_j$, gives
\[
\lambda\left(
\bigcup_{q\in\mathcal Q(A_j,B_j)}
\bigcup_{\gcd(a,q)=1}
\left\{
x\in[0,1]:
\left|x-\frac{a+\theta}{q}\right|
<
\frac{1}{8P(A_j,B_j)}
\right\}
\right)
<
2^{-B_j-1}.
\]
Since $B_1<B_2<\cdots$ are positive integers, we have $B_j\ge j$.
Consequently,
\[
\sum_{j=1}^{\infty}
2^{-B_j-1}
\le
\sum_{j=1}^{\infty}
2^{-j-1}
<
\infty.
\]
Thus
\[
\sum_{j=1}^{\infty}
\lambda\left(
\bigcup_{q\in\mathcal Q(A_j,B_j)}
\bigcup_{\gcd(a,q)=1}
\left\{
x\in[0,1]:
\left|x-\frac{a+\theta}{q}\right|
<
\frac{\psi(q)}q
\right\}
\right)
<
\infty.
\]
By the Borel--Cantelli lemma,
\[
\lambda\bigl(W_\theta^*(\psi)\bigr)=0.
\]

Thus, the inhomogeneous Duffin--Schaeffer conjecture fails for every
$\theta\in\mathcal G$. Since $\mathcal G$ is a dense $G_\delta$
subset of $\mathbb R$, the set of shifts for which the inhomogeneous
Duffin--Schaeffer conjecture fails is residual in $\mathbb R$.

\subsection*{Acknowledgements}
The counterexample was discovered with the assistance of OpenAI's GPT-5.6 Sol model; all mathematical arguments and references were verified by the authors. After completing the present manuscript, we shared it with Manuel Hauke-Treuer. We subsequently learned from him that Andrew Pollington had announced a  counterexample at a conference in York more than one year earlier.

Y. He was supported by the NSFC (No. 12401108) and partially by a grant from the Guangdong Provincial Department of Education (2025KCXTD013). We thank Manuel Hauke-Treuer for sharing their manuscript with us.


%

\end{document}